\documentclass[]{cJAS2eArXiv}

\usepackage{subfigure}% Support for small, `sub' figures and tables. Your choice of alternative may be preferred instead.

\theoremstyle{plain}

\usepackage{color}

\usepackage{amsmath}
\usepackage{mathrsfs}
\usepackage{bm}
\usepackage{graphicx}

\newcommand{\boma}[1]{\mbox{\boldmath ${#1}$}}

\begin{document}

%\jvol{00} \jnum{00} \jyear{2013} \jmonth{December}

%\articletype{GUIDE}

%\title{{\itshape Journal of Applied Statistics} \break \LaTeX\ style guide for authors (Style 2 + Reference Style S)}

\title{Interpretation of Compositional Regression with Application to Time Budget Analysis}

\author{
 I. M\"uller$^{\rm a}$,
  K. Hron$^{\rm a}$$^{\ast}$,
	\thanks{$^\ast$Corresponding author. Email: hronk@seznam.cz
\vspace{6pt}}
 E. Fi\v serov\'a$^{\rm a}$, 
 J. \v Smahaj$^{\rm b}$,
 P. Cakirpaloglu$^{\rm b}$,
 J. Van\v c\'akov\'a$^{\rm c}$\\
\vspace{6pt}  $^{\rm a}${\em Department of Mathematical Analysis and Applications of Mathematics, Faculty of Science, Palack\'y University, 17. listopadu 12, CZ-77146 Olomouc, Czech Republic};
$^{\rm b}${\em Department of Psychology, Philosophical Faculty, Palack\'y University, Vod\'arn\'i 6, 771 80 Olomouc, Czech Republic};
$^{\rm c}${\em Prostor Plus, Na Pustin\v e 1068, 280 02 Kol\'in, Czech Republic.}
}

\maketitle

\begin{abstract}

Regression with compositional response or covariates, or even regression between parts of a composition, is frequently employed in social sciences. Among other possible applications, it may help to reveal interesting features in time allocation analysis. As individual activities represent relative contributions to the total amount of time, statistical processing of raw data (frequently represented directly as proportions or percentages) using standard methods may lead to biased results. Specific geometrical features of time budget variables are captured by the logratio methodology of compositional data, whose aim is to build (preferably orthonormal) coordinates to be applied with popular statistical methods. The aim of this paper is to present recent tools of regression analysis within the logratio methodology and apply them to reveal potential relationships among psychometric indicators in a real-world data set. In particular, orthogonal logratio coordinates have been introduced to enhance the interpretability of coefficients in regression models.

\end{abstract}


\begin{classcode}
\textit{Classification codes}: 62J05, 62J12
\end{classcode}

\section{Introduction}

Regression analysis becomes challenging when compositional data as observations carrying relative information \cite{A86,PB11} occur in the role of response or explanatory variables. Although this might frequently seem to be a purely numerical problem, compositional data in any form inducing a constant sum constraint (proportions, percentages) rather represent a conceptual feature. In fact, compositional data may not necessarily be expressed with a constant sum of components (parts). The decision whether data at hand are compositional or not depends on the purpose of analysis - whether it is absolute values of components, or rather their relative structure, that is of primary interest.

One of most natural examples of compositional data are time budget (time allocation) data, discussed already in the seminal book on compositional data analysis \cite{A86}, p. 365. Apart from the compositional context, due to its psychological, social, and economic impacts, time allocation and its statistical analysis receives attention in many publications. The distribution of the total amount of time among productive-, maintenance-, and leisure activities reflects the current status and soundness of economy, with its labour-saving inventions, communication technologies, means of transportation, information and mass media channels, and level of consumption \cite{B65, J91, RG97, G00, G02}. The economy is usually closely linked to political arrangement, which through welfare state institutions (including child-care facilities) relieve citizens of many obligations, thus opening possibilities for loosening and restructuring their daily schedules \cite{K00, GS03, CL06}. Leisure time service is further provided for by various sports programs, holiday resorts, outdoor activities and the like, for both adolescents and adults.
Moreover, frequently also supplementary qualitative/quantitative variables (age, gender, variables resulting from psychometric scales) are of simultaneous interest, which calls for the use of regression modeling.

%Tradition, religion, and culture also play a prominent role in how people value different aspects of their lives, with direct repercussions for time use. Intriguing differences are brought forward by comparing the white Protestant U.S. \cite{S93} and Japanese \cite{F89} work-oriented values versus Latin America's focus on togetherness and friendship \cite{GL09}. On the side of an individual, socio-economic status, age, gender, and place of residence stand out among the influential factors. For example, Larson \& Verma \cite{LV99} and Larson et al. \cite{L01} give details on time budget for the age groups of children and young adolescents in different cultural settings. Moreover, \v S\'ipek \cite{S01} argues that free time will not be automatically experienced in a positive way. One has to be prepared to use their free time. One has to be able to admit their needs, to take a rest and accept a reward (like recreation, tourism, self-development).

When considering the problem of time allocation from the statistical point of view, the individual activities represent relative contributions to the overall time budget. Particularly, although the input data can be obtained either in the original time units, or directly in proportions or percentages, the relevant information is conveyed by ratios between the parts (time activities). Consequently, also differences between relative contributions of an activity should be considered in ratios instead of absolute differences as they better reflect relative scale of the original observations. 

Both scale invariance and relative scale issues are completely ignored when the raw time budget data or any representation thereof (like proportions or percentages) are analyzed using standard statistical methods. Although there do exist methods whose aim is to solve purely numerical problems resulting from the nature of observations carrying relative information (being of one dimension less than the actual number of their parts), these methods usually do not represent a conceptual solution to the problem of compositional data analysis. Instead, any reasonable statistical methodology for this kind of observations should be based on ratios between parts, or even \textit{logratios} (logarithm of ratios), which are mathematically much easier to handle \cite{A86, PB11}. Logratios as a special case of a more general concept of logconstrasts are used to construct coordinates with respect to the Aitchison geometry that captures all the above mentioned natural properties of compositions. Nevertheless, possibly due to apparent complexity of the logratio methodology, logratio methods haven't still convincingly entered psychometrical literature; methods to analyze time budget, mentioned in the seminal book of Van den Ark \cite{A99} and resulting from fixing the unit-sum constraint of compositional data, were mostly overcome during the last 15 years of intensive development in the field of compositional data. Very recently statistical analysis of psychological (ipsative) data seems to attract attention \cite{batista15,eijnatten15}. Nevertheless, still rather specific methods are used without providing a concise data analysis, particularly concerning regression modeling that frequently occurs in psychometrics.

For this reason, the aim of this paper is to perform a comprehensive regression analysis of time budget structure of college students by taking real-world data from a large psychological survey at Palack\'y University in Olomouc (Czech Republic). With that view, relations with other response/explanatory variables (as well as those within the original composition) will be analyzed using proper regression modeling.

The structure of the paper is as follows. In the next section, the orthonormal logratio coordinates are introduced first, and then regression modeling is discussed in more detail in Section 3. In order to achieve better interpretability of regression parameters while preserving all important features of regression models for compositional data, orthogonal coordinates (instead of orthonormal ones) are introduced as an alternative in Section 4. Section 5 is devoted to logratio analysis of the concrete time budget data set and the final Section 6 (Discussion) concludes.

\section{Orthonormal Logratio Coordinates for Compositional Data}

For a $D$-part composition $\mathbf{x}=(x_1,\ldots,x_D)'$, considering all possible logratios $\ln(x_i/x_j),\,i,j=1,\ldots,D,$ for statistical analysis means to take into account $D(D-1)/2$ variables (up to sign of the logarithm). This would lead to a complex ill-conditioned problem already for data sets with moderate number of variables. Moreover, information related to the original parts (although expressed possibly in logratios) is usually of primary interest. For this reason, a natural choice is to aggregate logratios meaningfully to logcontrasts (variables of type $\sum_{i=1}^Dc_i\ln x_i$, where $\sum_{i=1}^Dc_i=0$), that are able to capture all the relative information about single compositional parts (time activities). In other words, when $x_1$ plays the role of such a part, we proceed to variable $\ln(x_1/x_2)+\ldots +\ln(x_1/x_D)=(D-1)\ln(x_1/\sqrt[D-1]{\prod_{i=2}^Dx_i})$, i.e. to logcontrast that highlights the role of $x_1$ \cite{FH11}. In order to build a system of orthonormal coordinates, this variable needs to be further scaled and also the remaining $D-2$ coordinates, orthonormal log-contrasts, are constructed consequently (we refer to isometric logratio (ilr) coordinates \cite{E03}). One possible choice of ilr coordinates that fulfill the above requirements (for any of parts $x_l,\,l=1,\dots,D,$ in place of $x_1$) is $\mathbf{z}^{(l)}=(z_1^{(l)},\dots,z_{D-1}^{(l)})'$, 
\begin{equation}
\label{eilr}
z_i^{(l)}=\sqrt{\frac{D-i}{D-i+1}}\,\ln
\frac{x_{i}^{(l)}}{\sqrt[D-i]{\prod_{j=i+1}^Dx_j^{(l)}}}, 
\ i=1,\ldots ,D-1 .
\end{equation}
The case of $x_1$ would be obtained by choosing $l=1$. In a more general setting, the composition $(x_1^{(l)},x_2^{(l)},\ldots,x_l^{(l)},x_{l+1}^{(l)},\dots,x_D^{(l)})'$ stands for such a permutation of the parts $(x_1,\dots,x_D)'$ that always the $l$-th compositional part fills the first position, $(x_l,x_1,\dots,x_{l-1},x_{l+1},\dots,x_D)'$. In such a configuration, the first ilr coordinate $z_1^{(l)}$ explains all the relative information (merged into the corresponding logcontrast) about the original compositional part $x_l$, the coordinates $z_2^{(l)},\dots,z_{D-1}^{(l)}$ then explain the remaining logratios in the composition. Note that the only important position is that of $x_1^{(l)}$ (because it can be fully explained by $z_1^{(l)}$), the other
parts can be chosen arbitrarily, because different ilr coordinates are orthogonal rotations of each other \cite{E03}. Although this particular choice of ilr coordinates has been used successfully in many geological and chemometrical applications \cite{BEP14, F12, K14}, no experiences are recorded in the psychometrical context.

\section{Regression Analysis within the Logratio Methodology}
 
% Regression 
Regression analysis is an important tool for analyzing the relationships between the response variable $Y$ and known explanatory variables $\mathbf{x}$, see, e.g. \cite{MPV06}. Although in the psychometrical context it is often difficult to distinguish whether the covariates are driven by an error as well, or not, we will follow the assumption of fixed covariates in order to enable estimation of regression parameters using the standard least squares (LS) method, resulting in easy-to-handle statistical inference (hypotheses testing). When the response variables or explanatory variables are compositional, special treatment in regression is necessary. A natural way for introducing regression with compositional explanatory variables $\mathbf{x}=(x_{1},x_{2},\ldots,x_{D})'$  is to perform a standard multiple regression where the explanatory variables $\mathbf{z}_i=(1,z_{i,1},z_{i,2},\ldots,z_{i,D-1})'$ represent the ilr coordinates of $\mathbf{x}_i$ and 1 for the intercept \cite{H12}. Using a special choice  of ilr coordinates $\mathbf{z}^{(l)}$ given by (\ref{eilr}), we can consider the $l$th ilr basis, for $l=1,2,\ldots,D$, and we obtain $D$ different multiple regression models in the form 
\begin{equation}
\label{Xcomp}
   Y_i=\beta_0+\beta_1^{(l)} z_{i,1}^{(l)}+\cdots+\beta_{D-1}^{(l)}z_{i,D-1}^{(l)}+\varepsilon_i^{(l)},\ i=1,2,\ldots,n,
\end{equation}
where $\beta_0, \beta_1^{(l)},\ldots \beta_{D-1}^{(l)}$ are unknown regression parameters and $\varepsilon_i^{(l)}$ are random errors in the $l$th model. Due to the orthogonality of different ilr bases, the intercept term $\beta_0$ is the same for all $D$ models (similarly as the index of determination $R^2$ or the $F$ statistic to test the overall significance of the covariates) \cite{H12}. The regression parameters can be estimated in the standard way by the least squares (LS) method. Using the notation $\mathbf{Y}=(Y_1,Y_2,\ldots,Y_n)'$ for the observation vector, $\mathbf{Z}^{(l)}=(\mathbf{z}_1^{(l)}, \mathbf{z}_2^{(l)},\ldots,\mathbf{z}_n^{(l)})'$ for $n\times D$ design matrix, $\boma{\beta}^{(l)}=(\beta_0,\beta_1^{(l)},\ldots,\beta_{D-1}^{(l)})'$
for regression parameters,  and $\boma{\varepsilon}^{(l)}=(\varepsilon_1^{(l)},\varepsilon_2^{(l)},\ldots,\varepsilon_n^{(l)})'$ for the error term, models (\ref{Xcomp}) can be rewritten in the matrix form 
\begin{equation}
\label{Xcomp2}
   \mathbf{Y}=\mathbf{Z}^{(l)}\boma{\beta}^{(l)}+\boma{\varepsilon}^{(l)},\ l=1,2,\ldots,D.
\end{equation}
We can consider that random errors in the $l$th model are not correlated with the same variance $\sigma^2_{(l)}$. Then
the best linear unbiased estimators of regression parameters $\boma{\beta}^{(l)}$ by the LS method are
\begin{equation}
   \widehat{\boma{\beta}}^{(l)}=
		   (\mathbf{Z}'^{(l)}\mathbf{Z}^{(l)})^{-1}\mathbf{Z}'^{(l)} \mathbf{Y}, 
		    \ l=1,2,\ldots,D.
\end{equation}  
From the practical point of view, only the parameter $\beta_1^{(l)}$ is important, since it corresponds to the first ilr coordinate $z_{1}^{(l)}$ that explains all the relative information about the part $x_1^{(l)}$. The other parameters $\beta_2^{(l)},\ldots,\beta_{D-1}^{(l)}$ do not have such straightforward interpretation. So, we can say, e.g., that the absolute change of the conditional mean of $Y$ with respect to coordinate $z_1^{(l)}$ is about $\beta_1^{(l)}$, if other coordinates $z_j^{(l)}$, $j=2,3,\ldots,D-1$ (representing subcomposition $(x_1,\dots,x_{l-1},x_{l+1},\dots,x_D)'$), are fixed.

The unbiased estimator of $\sigma^2_{(l)}$ in the $l$th model (\ref{Xcomp2}) is
\begin{equation}
    \widehat{\sigma}^2_{(l)}=
		   (\mathbf{Y}-\mathbf{Z}^{(l)}\widehat{\boma{\beta}}^{(l)})'(\mathbf{Y}-\mathbf{Z}^{(l)}\widehat{\boma{\beta}}^{(l)})/(n-D),
\end{equation}  
that can be used to estimate the variance-covariance matrix of the estimator of regression parameters,
\begin{equation}
\widehat{\mathrm{var}}(\widehat{\boma{\beta}}^{(l)})
    =\widehat{\sigma}^2_{(l)} (\mathbf{Z}'^{(l)}\mathbf{Z}^{(l)})^{-1}.
\end{equation}		
 
Under assumption of normality of random errors we can perform any standard statistical inference, e.g. test the significance of regression parameters, or to construct confidence intervals for them. The significance of the individual regression parameters in the $l$th model, $l=1,2,\ldots,D$, can be tested by the following statistics:
\begin{equation}
	T_0=\frac{\widehat{\beta}_0}{\widehat{\sigma}_{(l)}\sqrt{\{(\mathbf{Z}'^{(l)}\mathbf{Z}^{(l)})^{-1}\}_{1,1}}};\quad
	T_i^{(l)}=\frac{\widehat{\beta}_i^{(l)}}{\widehat{\sigma}_{(l)}\sqrt{\{(\mathbf{Z}'^{(l)}\mathbf{Z}^{(l)})^{-1}\}_{i+1,i+1}}},\quad
	i=1,2,\ldots,D-1.
\end{equation}
Here the symbol $\{(\mathbf{Z}'^{(l)}\mathbf{Z}^{(l)})^{-1}\}_{i+1,i+1}$ denotes the $(i+1)$th diagonal element of the matrix $(\mathbf{Z}'^{(l)}\mathbf{Z}^{(l)})^{-1}$. Under the null hypothesis that regression parameters are zeros, the statistics $T_0$ and $T_i^{(l)}$ each follow a Student $t$-distribution with $n-D$ degrees of freedom. The statistic $T_0$ is the same irrespective of the choice of $l=1,\dots,D$ in (\ref{Xcomp}), see \cite{H12} for details. Of course, the response variable can have also another distribution than normal, i.e. the methodology of generalized linear models \cite{DB08} can be directly implemented.

Similarly, when the response variables $\mathbf{Y}=(Y_{1},Y_{2},\ldots,Y_{D})'$ are compositional and explanatory variables $\mathbf{x}=(x_{1},x_{2},\ldots,x_{k})'$ are non-compositional, one can use the regression models where the response variables $z_{j}$, $j=1,2,\ldots,D-1$,  represent the ilr coordinates of $\mathbf{Y}$ \cite{E11}. Using the ilr coordinates (\ref{eilr}), where only the first ilr coordinate $z_1^{(l)}$ is of interest, we obtain $D$ different multiple regression models in the form 
\begin{equation}
\label{Ycomp}
   Z_{i,1}^{(l)}=\gamma_0^{(l)}+x_{i,1}\gamma_1^{(l)} +\cdots+x_{i,k}\gamma_{k}^{(l)}+\varepsilon_i^{(l)},\ i=1,2,\ldots,n,\ l=1,2,\ldots,D.
\end{equation}
In this case, the interpretation of regression parameters is the following. For example, if $x_2,\ldots,x_k$ are fixed, then for each change of 1 unit in $x_1$, the conditional mean of $Z_1^{(l)}$ changes $\gamma_1^{(l)}$ units. Nevertheless, similarly as for the case of regression with compositional explanatory variables, because the orthonormal coordinates (\ref{eilr}) have to be interpreted in terms of \textit{scaled} logratios under natural logarithm, the interpretation of these ``units'' and thus also values of regression parameters might get rather complex for practical purposes. Under the usual multiple regression model assumptions, (\ref{Ycomp}) can be expressed in the matrix form
\begin{equation}
\label{Ycomp2}
   \mathbf{Z}_1^{(l)}=\mathbf{X}\boma{\gamma}^{(l)}+\boma{\varepsilon}^{(l)},\ l=1,2,\ldots,D,
\end{equation}
where $\mathbf{Z}_1^{(l)}=(Z_{1,1}^{(l)},Z_{2,1}^{(l)},\ldots,Z_{n,1}^{(l)})'$ is an observation vector, $\mathbf{\gamma}=(\gamma_0,\gamma_1,\ldots,\gamma_k)'$ is a vector of regression parameters, and $\mathbf{X}=(\mathbf{1},\mathbf{x}_1, \mathbf{x}_2,\ldots,\mathbf{x}_k)$ is $n\times (k+1)$ design matrix. Here $\mathbf{1}$ is a vector of $n$ ones. When the random errors in the $l$th model are not correlated with the same variance $\sigma^2_{e,(l)}$, the best linear unbiased estimator of regression parameters $\boma{\gamma}^{(l)}$ by the LS method is
\begin{equation}
   \widehat{\boma{\gamma}}^{(l)}=
		   (\mathbf{X}'\mathbf{X})^{-1}\mathbf{X}\mathbf{Z}_1^{(l)},
		    \ l=1,2,\ldots,D,
\end{equation} 
with the estimated variance-covariance matrix
\begin{equation}
\widehat{\mathrm{var}}(\widehat{\boma{\gamma}}^{(l)})
    =\widehat{\sigma}^2_{e,(l)} (\mathbf{X}'\mathbf{X})^{-1}.
\end{equation}
The unbiased estimator of $\sigma^2_{e,(l)}$ in model (\ref{Ycomp2}) is 
\begin{equation}
    \widehat{\sigma}^2_{e,(l)}=
		   (\mathbf{Z}_1^{(l)}-\mathbf{X}\widehat{\boma{\gamma}}^{(l)})'(\mathbf{Z}_1^{(l)}-\mathbf{X}\widehat{\boma{\gamma}}^{(l)})/(n-k-1).
\end{equation}

Again, under assumption of normality of random errors we can test the significance of regression parameters, or construct confidence intervals for them. In this case, the significance of the individual regression parameters in the $l$th model, $l=1,2,\ldots,D$, can be tested by the  statistic:
\begin{equation}
	U_i^{(l)}=\frac{\widehat{\gamma}_i^{(l)}}{\widehat{\sigma}_{e,(l)}\sqrt{\{(\mathbf{X}'\mathbf{X})^{-1}\}_{i+1,i+1}}},\quad
	i=0,1,\ldots,k.
\end{equation}
Under the null hypothesis that regression parameters are zeros, the statistics $U_i^{(l)}$ follow a Student $t$-distribution with $n-k-1$ degrees of freedom. The confidence intervals for regression parameters on confidence level $1-\alpha$ in the $l$th model, $l=1,2,\ldots,D,$ would be constructed analogously as before.
% are
%\begin{equation}
%	\mathcal{I}_{1-\alpha}(\gamma_i^{(l)})
%	       =\left\langle 
%	            \widehat{\gamma}_i^{(l)}\pm\widehat{\sigma}_{e,(l)}\sqrt{\{(\mathbf{X}'\mathbf{X})^{-1}\}_{i+1,i+1}}t_{n-k-1}(1-\alpha/2)	
%				  \right\rangle,\
%	i=0,1,\ldots,k.							
%\end{equation}

Finally, within the logratio methodology we can consider also the case of regression among parts of a composition, in particular, between a part $x_0$ and the rest of compositional parts, $x_1,\dots,x_D$, in a $(D+1)$-part composition. Following \cite{BEP14, HTHF15}, a natural choice is to consider the case of regression with compositional explanatory variables, where the response is formed by coordinate, carrying the relative information of $x_0$ (with respect to compositional covariates), i.e.,
$$z_0=\sqrt{\frac{D}{D+1}}\ln\frac{x_0}{\sqrt[D]{\prod_{i=1}^Dx_i}}.$$
By construction, $z_0$ is orthonormal to the rest of coordinates, assigned to explanatory parts as in (\ref{eilr}).

\section{Orthogonal Coordinates for Compositional Regression}
\label{orthogcoor}

Although the above regression models in orthonormal logratio coordinates are theoretically well justified, both the normalizing constants to reach orthonormality and the natural logarithm itself result in quite a complex interpretation of the regression parameters. A way out is to move to \textit{orthogonal} coordinates, where nothing from the above properties of regression modeling in coordinates is lost (in particular, values of $T_i^{(l)}$, $U_i^{(l)}$ and $F$ statistics, neither the geometrical features of regression with compositional response \cite{E11}), while, at the same time, a substantial simplification in parameter interpretation is gained. Following (\ref{eilr}), these considerations lead to orthogonal coordinates
\begin{equation}
\label{eilrog}
z_i^{(l)*}=\log_2\frac{x_{i}^{(l)}}{\sqrt[D-i]{\prod_{j=i+1}^Dx_j^{(l)}}}, 
\ i=1,\ldots ,D-1 ,
\end{equation}
for $l=1,\dots, D$, where the normalizing constants are omitted and the original natural logarithm is replaced by the binary one. Let's see the effect of using the orthogonal coordinates for all regression models introduced above (parameters of their corresponding versions in orthogonal coordinates (\ref{eilrog}) are always marked with an asterisk). Considering regression with compositional explanatory variables first, from properties of LS estimation and the relation between logarithms of different bases we get 
$$\beta_0^*=\beta_0,\ \beta_1^{(l)*}=\ln(2)\sqrt{\frac{D-1}{D}}\beta_1^{(l)},$$
generally
$$\beta_i^{(l)*}=\ln(2)\sqrt{\frac{D-i}{D-i+1}}\beta_i^{(l)},\ i=1,\ldots,D-1,$$
and similarly for their estimates and the respective standard errors. Analogously, for models resulting from regression with compositional response we get
$$\gamma_i^{(l)*}=\log_2(\mathrm{e})\sqrt{\frac{D}{D-1}}\gamma_i^{(l)},\ i=0,\ldots, k.$$
Finally, in regression within composition both the above effects are combined, i.e., for $D$ regression models 
\begin{equation}
\label{YXcomp}
   Z_{i0}=\beta_0+\beta_1^{(l)} z_{i,1}^{(l)}+\cdots+\beta_{D-1}^{(l)}z_{i,D-1}^{(l)}+\varepsilon_i^{(l)},\ i=1,2,\ldots,n,
\end{equation}
($l=1,\dots,D$) we obtain
$$\beta_0^*=\log_2(\mathrm{e})\sqrt{\frac{D+1}{D}}\beta_0,\ \beta_i^{(l)*}=\sqrt{\frac{(D+1)(D-i)}{D(D-i+1)}}\beta_i^{(l)},\ i=1,\ldots,D-1.$$

Indeed, the interpretation of regression coefficients gets simpler now. 
For regression with compositional regressors and non-compositional response, first note that a unit additive increment in a log-transformed coordinate $z$ is equivalent to a two-fold multiplicative increase in the relative dominance of the original compositional variable $x$, if the base-2 logarithm is used, that is,
$$\Delta z_1^{(l)*}=\log_2\frac{x_1^{(l)*}}{\sqrt[D-1]{\prod_{i=2}^Dx_i^{(l)*}}}\cdot2-\log_2\frac{x_1^{(l)*}}{\sqrt[D-1]{\prod_{i=2}^Dx_i^{(l)*}}} = 1.$$
The coefficient $\beta_1^{(l)*}$ in the regression equation then has the usual meaning of an additive increase in the response $y$
that corresponds to increasing $z$ by one (i.e., increasing the dominance of $x$ twice), while keeping all else fixed. For example, if $\beta_1^{(l)*}=3$, the value of the response gets higher by 3 units when the relative dominance of the part $x_l$ with respect to the average of the other parts, see the logratio in (\ref{eilrog}), is doubled, at constant values of the other involved covariates (orthogonal coordinates). Next, in case of regression with compositional response and non-compositional regressors, $\gamma_j^{(l)*}$ is the additive increment of the log-transformed response $z$ when adding one to an explanatory variable $x_j,\,j=1,\dots,k,$ (at constant values of the other covariates)
$$\gamma_j^{(l)*} = \Delta Z_1^{(l)*}=\log_2\frac{Y_1^{(l)*}}{\sqrt[D-1]{\prod_{i=2}^DY_i^{(l)*}}}\delta-\log_2\frac{Y_1^{(l)*}}{\sqrt[D-1]{\prod_{i=2}^DY_i^{(l)*}}}=\log_2\delta,$$
where $\delta=2^{\gamma_j^{(l)*}}$  is the multiplicative increase in the relative dominance of the original compositional response $y$. So, for a unit additive change in $x_j$, the ratio of $Y_1^{(l)}$ to the ``mean value'' of the other compositional responses grows $\delta=2^{\gamma_j^{(l)*}}$ times. Finally, an analogous interpretation for regression within composition can be obtained, namely, a two-fold multiplicative increase in the relative dominance of $x_l$ (or equivalently, a unit additive increment in coordinate $z_1^{(l)*}$) brings about the increase in the relative dominance of response $x_0$ of 
$$\delta=2^{\beta_1^{(l)*}},\mbox{ where }\Delta z_0^*=\log_2\delta.$$
Note also that the above expression for the proportionality coefficient $\delta$ stays the same irrespective of the base to which the logarithm was taken, as factor 2 in the expression now stands for two-fold increase in dominance, not for the logarithmic base.

\section{Time Budget Analysis}

Following the previous developments, the decision to admit that the time budget data are by their nature compositional invites one to couch analysis in terms of logratios instead of working with the original observations in percentages; namely, the latter would lead to biased conclusions due to relative character of compositions. The aim of this section is to demonstrate on real-world psychometric data that working with logratios in the 
regression context is as accessible as dealing with the original observations.

\subsection{Data and Methods}

For this purpose, we employ data from \cite{V13} that were obtained in a large psychometric study, guaranteed and realized by the Department of Psychology, Palack\'y University in Olomouc, Czech Republic. A questionnaire called ``Leisure Time" was distributed among students at the above university, reaching a total of $N$ = 414 respondents (347 women, 67 men) who provided complete answers. The items included in the questionnaire tapped three distinct areas: i) personal characteristics (age, gender, faculty and field of study); ii) leisure time (its concept, absolute and relative amount, content); iii) personality traits (self-esteem and attitude to challenges). In terms of current analysis, of particular interest are relationships among the following variables: Daily Time Budget as expressed in seven compositional variables (parts, summing up to 100 percent) \textit{study/work}, \textit{commuting}, \textit{food}, \textit{hygiene\&dressing}, \textit{sleep}, \textit{household duties}, and \textit{leisure time}; personality variables \textit{self-esteem} (z-score from a 10-item Rosenberg Self-Esteem Scale \cite{R65} included in the questionnaire) and \textit{challenge} (``Are you a person who invites challenges, i.e. opportunities to surpass yourself?", originally 4-choice response collapsed into dichotomic and coded as 1 for ``always" or ``almost always", and 0 for ``almost never" or ``never"); and covariates of \textit{age} (in years) and \textit{gender} (dichotomic, coded as 1 for men and -1 for women). Distribution of the variables \textit{age} and \textit{self-esteem} is visualized in Figure 1 in the form of EDA-plots using the R package StatDA \cite{F13}.

\begin{center}
\textit{Figure 1 near here.}
\end{center}

Although the respondents were asked to enter data on Daily Time Budget in percentages, the obtained range of the sum of parts was $\langle 7, 520\rangle$ due to misunderstanding the units to use and their prescribed constant sum constraint (of course, most of the row sums were exactly or close to 100). Nevertheless, the important information on relative contributions of parts to the overall time budget was unaffected by using whatever units, which thus emphasizes even further the necessity to apply the logratio methodology in statistical processing. Note once again that for the logratio methodology the constant sum representation of compositional data is not a necessary requirement. However, for the purpose of easier comparisons, in the following the percentage representation was taken for all time budget observations.
%percentage representation was adopted for all time budget observations as well. We refer to time budget data as compositional. 

Besides paying attention to differences, as well as agreement, in logratio vs. ``standard" methodologies, we will keep our thoughts focused on some tentative conjectures about interconnections among variables. This data set allows for exploring possible influences among several prominent psychological factors. On the one hand, we have the pair of personality traits of self-esteem and openness to challenge which we expect to be bundled close together and even boost each other if challenges are being tackled successfully, or else restrain each other in a downward spiral. On the other hand, the necessity of time allocation brings about an inevitable interplay of work, active relaxation, and sleep (passive relaxation). And then, of course, these two broad areas come into mutual contact in complex ways. 

These considerations lead us, at the outset, to postulate a firm and positive relationship between personality traits of \textit{challenge} and \textit{self-esteem}. Next, within compositional variables, we deem as highly probable a negative relationship between \textit{work/study} and \textit{leisure time}, and between \textit{work/study} and \textit{sleep} on the premise that working/studying takes away time from both these forms of relaxation. \textit{Sleep} is considered loosely associated with \textit{leisure time} on the grounds that the time left after deducting all duties is being distributed between both. If there is more time available, it will add up to both sleep and leisure. If any at all, the relationship between \textit{sleep} and \textit{challenge} is expected to be negative, as the person who is busy taking challenges might have less time for sleep. The association between \textit{sleep} and \textit{self-esteem} is less clear-cut but it can be conceived along the lines that a self-assured person participates in numerous activities and thus sleeps less, while, on the other hand, an insecure person may seek sleep as a welcome escape from reality. As a consequence, \textit{work/study} should be positively related with both \textit{challenge} and \textit{self-esteem}, and \textit{leisure time} negatively related with both. Any effects of gender may be obscured in this dataset as men are seriously underrepresented among respondents.

In the following, the relationships among variables are determined through regression analysis. A logratio approach (which is deemed appropriate whenever a compositional variable out of Daily Time Budget is present) is compared to a standard non-compositional approach, e.g. Linear Model (LM) or Generalized Linear Model (GLM). In the statistical analysis we focused on those relations that are primarily not gender related. Moreover, preliminary exploratory analysis using variation matrix \cite{A86} and compositional biplot \cite{AG02} revealed strong relationship between \textit{food} and \textit{hygiene\&dressing components; because of} their rather marginal importance for psychological interpretation, these parts will be excluded from further consideration (but kept as parts of the initial composition).

\subsection{Regression Analysis}

From the essence of the data set, interconnections among variables (compositional and non-compositional, or even within the time budget composition) are of primary interest. For this purpose, several regression models were applied to data. Accordingly, in addition to Daily Time Budget, non-compositional variables of \textit{challenge}, \textit{self-esteem}, \textit{age}, and \textit{gender} were taken into consideration here. In order to enable direct interpretation of regression output, orthogonal coordinates (as described in Section 4), instead of orthonormal ones, were employed for the compositional variables within logratio approach.

As a first step, let us explore the manner how seeking challenges is determined by Daily Time Budget and other explanatory variables. That is, the response now is non-compositional (binomial), while some of the regressors are compositional and others not. For this purpose, binomial regression (a special case of logistic regression) was applied, first with compositional regressors in logratio coordinates, second with the original variables in percentages; note that any representation of the orthogonal logratio coordinates would lead to the same parameter estimates for the non-compositional covariates. From the time budget variables just those of potential psychological influence were included (\textit{study/work}, \textit{commuting}, \textit{sleep}, \textit{household duties}, and \textit{leisure time}); of course, due to construction of the regression model in coordinates, all parts of the original composition were taken into account for the estimation purposes under logratio approach. On the other hand, perfect collinearity among compositional variables makes it impossible to include all of them simultaneously as regressors in a standard linear model. Following \cite{H12}, common regression output like parameter estimates, their standard errors, values of corresponding statistics and their P-values (using function \texttt{glm} from R-package MASS, see \cite{VR02} for further details) are collected in Table~\ref{tab:1} (all tables with detailed results are included as supplementary material), where names of the original parts stand as notation for the corresponding orthogonal coordinates (\ref{eilrog}). It can be seen that both the \textit{study/work} coordinate and the \textit{self-esteem} variable are contributing the most (in the positive direction, due to positive sign of their coefficients) in explaining the \textit{challenge} response. The interpretation of coefficients is such that if the relative dominance of \textit{study/work} in time budget doubles (with respect to average contribution of the other parts), there is a shift of .42 towards seeking challenges in the response binomial logit (other covariates staying fixed); and this effect is about the same for a unit increase in \textit{self-esteem} z-score (a shift of .45). Note that, in line with the methodology described in the previous section, five regression models were employed to obtain the estimates for the compositional coordinates. By applying orthogonal coordinates (\ref{eilrog}), the interpretation of regression coefficients gets much easier than with original orthonormal coordinates (\ref{eilr}). The tight link between \textit{challenge} and \textit{self-esteem} is thus established. On the other hand, we don't see significance of either \textit{sleep} or \textit{leisure time}, though the direction (sign of coefficient) is as expected.

\begin{center}
\textit{Table 1 near here.}
\end{center}

The output of binomial regression with the original compositional variables is shown in Table~\ref{tab:2}. In general, there is not much difference from the logratio approach above (also the model fit, expressed by AIC criterion, stays almost the same), which would indicate that the distortion of covariance structure among percentage covariates (see, e.g., \cite{A86} for details) didn't have dramatic influence on regression output. The strength of association between openness to \textit{challenge} and \textit{self-esteem} remains unchanged. Nevertheless, the interesting influence of \textit{study/work} coordinate, which was clearly visible using the logratio coordinates, is now lost.

\begin{center}
\textit{Table 2 near here.}
\end{center}

As a second step, let us look the other way around and search for possible significant covariates of Daily Time Budget. For this purpose regression with compositional response was employed, the response variables now being the five chosen Daily Time Budget variables. The logratio approach leads to five univariate regression models (with orthogonal coordinates corresponding to individual compositional parts) and the results are displayed in Table~\ref{tab:3} (to save space, just regression estimates and possible significance at the usual level $\alpha=0.05$, marked by asterisk, are provided). The effects of particular covariates on response coordinates are evident. For example, by increasing the value of \textit{self-esteem} by one, the relative dominance of \textit{leisure time} in the composition (with respect to average of parts) increases approximately by 6 percent ($2^{0.088} = 1.06$). Similarly, taking challenges brings the relative dominance of \textit{study/work} 18 percent higher ($2^{0.237} = 1.18$), and one more year of age 2.9 percent higher. The positive association between \textit{study/work} and taking \textit{challenges} is in accordance with our anticipations, but with \textit{self-esteem} and \textit{leisure time} a contrary direction was expected. The connection between \textit{sleep} and both \textit{challenge} and \textit{self-esteem} remained below significance. It is interesting to see also some gender influence on both \textit{sleep} and \textit{leisure time}. Due to coding used (1 for male and $-1$ for female) it can be concluded that for males \textit{sleep} and \textit{leisure time} play a more important role in the overall time budget than for females. More precisely, the part \textit{sleep} is explained only by gender. Hence, $\widehat{z}_1^{(sleep)*}=1.644$ is the fitted value of the coordinate ${z}_1^{(sleep)*}$ for males, while $\widehat{z}_1^{(sleep)*}=1.357$ for females. It means that the relative dominance of \textit{sleep} in the composition to the
``mean value'' of the other compositional responses is $2^{1.644}=3.125$ for males (3.125 times higher relative contribution of sleep than for the averaged rest of components), and $2^{1.357}=2.562$ for females. Further, it can be concluded that the relative dominance of \textit{sleep} for males is $2^{2\widehat{\gamma}_{(sleep)}^*}=1.22$ times greater than for females. Although results for \textit{food} and \textit{hygiene\&dressing} variables are in general not discussed in this section, it is worth to note that for \textit{hygiene\&dressing} a significant role of gender (in the negative sense) was revealed; accordingly, this compositional part plays a more important role in time budget of females than for males.

\begin{center}
\textit{Table 3 near here.}
\end{center}

By way of comparison, the same regression model was analyzed under the assumption of Dirichlet distribution for the compositional response that is popular also in psychometric context \cite{G08} and, although rather inconsistent with logratio methodology, is still frequently recommended for modeling compositional data. For this purpose function \texttt{DirichReg} from the package DirichletReg \cite{M14} was applied by expressing the input compositions in proportional representation; regression output is collected in Table~\ref{tab:4}. Apart from its apparent computational complexity, Dirichlet regression does not seem to shed new light into the problem; moreover, some of the potential relationships that have emerged with the logratio approach are lost again.

\begin{center}
\textit{Table 4 near here.}
\end{center}

From the previous analysis, \textit{leisure time} seems to be strongly linked with the non-compositional variables. A natural question thus arises whether regression could reveal also some relations within parts of the time budget composition. Thus, as the third step, the corresponding logratio model from Section 2 was applied, by expressing both the response and regressors in orthogonal logratio coordinates (and with additional non-compositional covariates). Similarly as before, Table~\ref{tab:5}  collects results from four regression models, each highlighting the role of one of compositional explanatory variables (without influence on the non-compositional covariates). Though the $R^2$ statistic gives rather low value (as is usual in social science), some patterns stand out. In particular, relative dominance of \textit{leisure time} is positively influenced by \textit{sleep} (increasing the dominance of sleep twice enlarges the dominance of leisure time by 27 percent, as $2^{0.34} = 1.27$) and marginally by \textit{self-esteem} (double increase in self-esteem increases the dominance of leisure time by 6 percent); negative effects on leisure time are formed by \textit{study/work} (10 percent decrease in dominance, $2^{-0.16} = 0.90$) and \textit{commuting} (decrease in dominance by 13 percent). Consistency with the previous logratio model (Table~\ref{tab:3}, regression with compositional response) is underlined by the roles of \textit{self-esteem} and \textit{gender} covariates. Again, a psychological interpretation can be easily derived. Here we are able to pinpoint the significant positive association of \textit{sleep} and \textit{leisure time}, as well as negative association of \textit{work/study} and \textit{leisure time}. Marginally significant is the connection between \textit{leisure time} and \textit{self-esteem} which appeared significant in previous regression (Table~\ref{tab:3}).

\begin{center}
\textit{Table 5 near here.}
\end{center}

For the final comparison we consider the standard linear regression model where the original parts in percentages are involved (except for \textit{food} and \textit{hygiene\&dressing}), see Table~\ref{tab:6}  for the regression summary. Although conclusions from this model as regards non-compositional covariates would be pretty similar as with logratio methodology, the situation is different in other respects. By comparing $R^2$ for these two models and P-values at respective compositional covariates it is easy to see that for the standard regression model these values are very strongly driven by the constant-sum constraint of the original composition. In particular, note that by including all the compositional parts, $R^2$ would be brought up to 1, i.e., relations between the response and covariates would be completely driven by constant sum constraint of the input data. Of course, as statistical processing of the original compositions violates both scale invariance and relative scale properties of observations, it cannot be concluded that by considering compositional data without a constant sum constraint, the resulting regression model would be relevant. Nevertheless, in percentage representation, which is the case here, the irrelevance of the standard approach is clearly observable.

\begin{center}
\textit{Table 6 near here.}
\end{center}

\subsection{Results}

The logratio approach to regression analysis supports our hypothesis of strong negative association between \textit{work/study} and \textit{leisure time}, as well as of strong positive association between \textit{challenge} and \textit{self-esteem}. Next, \textit{leisure time} is significantly tied to  \textit{self-esteem} but the direction here appeared positive, rather than negative as expected. The reason could be that self-assured people don't feel the urge to work that much and rather take things easy, allowing themselves more leisure. Also, an explanation in keeping with \cite{S01} says that people with higher self-esteem may be better prepared to use their free time and it may be easier for them to admit their needs (for rest and reward). The connection between \textit{leisure time} and \textit{challenge} was not born out (remained below significance, though direction was negative as anticipated). The above regression results were agreed on by both logratio and standard linear model approaches. Both approaches also showed a relationship between \textit{sleep} and \textit{leisure time}. However, here the directions differed: logratio showed it to be positive (as hypothesized), linear model negative. On top of that, logratio approach was capable of revealing a significant positive connection between \textit{challenge} and \textit{work/study}.

The psychologically relevant variables seem to form a well-defined cluster of \textit{challenge, self-esteem}, and \textit{work/study}. Somewhat in opposition stands the pair of \textit{leisure time} and \textit{sleep}. However, their position with respect to the main cluster is less clearly marked, as \textit{leisure time} is negatively linked to \textit{work/study} but positively (perhaps only marginally) to \textit{self-esteem}. Nevertheless it seems reasonable to assume that working/studying does take time away from both leisure and sleep simultaneously. 

Finally, it is also worth noting that standard regression models were presented mostly for the sake of comparison of the logratio approach with alternatives that would be most possibly used instead. While in some cases their output might seem meaningful, it can also happen that by ignoring the relative structure of Daily Time Budget some interesting features are lost, as was the case in Table \ref{tab:2} and Table \ref{tab:4}. For some cases, like when percentage representation of the relative contributions is analyzed, it is very easy to demonstrate that scale invariance of compositional data leads to clear failure of the standard approach (Table \ref{tab:6}).

\section{Discussion}

Specific habits of time allocation reveal a lot about an individual, a community, a society, or a culture. In each society, options available to individuals for earning their living determine the amount of time they will spend working, or preparing themselves for any such productive activity through study or apprenticeship. In modern times, we have witnessed a continuous reduction in working hours, at least in industrialized countries. At the same time, due to constant total time budget, this development leaves more space for other activities, both necessary (self- and home-maintenance like sleep, eating, hygiene, care for family and house) and discretionary (leisure activities like socializing, culture, sports, reading, idling, etc.). As the time budget data are usually accompanied with other psychometric variables, regression modeling is the first and intuitive choice for a relevant statistical analysis.

Due to relative character of time budget allocation, it seems natural to work with (log-)ratios rather than with observations in the original scale (i.e. represented usually in proportions or percentages). It turned out that logratios meet the scale invariance and relative scale requirements (among others that are important for reasonable processing of compositional data) commonly raised in connection with any observations carrying primarily relative information. The main problem is then how to construct logratio coordinates, both meaningful from the mathematical point of view (guaranteed in particular by orthonormality of coordinates) and at the same time providing easy interpretation. The aim of the paper was to enhance interpretability of regression analysis output by employing orthogonal coordinates in place of the mathematically preferred orthonormal ones, demonstrated for the particular case of time budget data. The reason for the choice of alternative coordinates is that all the beneficial properties of the orthonormal coordinates are maintained also by the orthogonal ones, but the latter enable (by avoiding the scaling constants and changing the logarithmic base) a more straightforward interpretation. We are convinced that better interpretability of the regression models, discussed in the paper, can help with applicability of the logratio methodology in psychological research, and also in general.

\label{lastpage}

\def\refname{References}

\vspace{\fill}\newpage

\begin{table}[htp]
%\label{tab:cancerlm2}
\caption{\label{tab:1}Logratio approach: Results from regression of \textit{challenge} on orthogonal coordinates of the explanatory composition and further covariates. For explanations see text.}
\begin{center}
\begin{footnotesize}
\begin{minipage}{10cm}
\begin{verbatim}

                   Estimate Std. Error z value Pr(>|z|)
(Intercept)       -0.69708    1.24154  -0.561  0.57448
study/work         0.42200    0.14164   2.979  0.00289
commuting         -0.06723    0.10961  -0.613  0.53959
sleep             -0.20460    0.17476  -1.171  0.24168
household duties  -0.02904    0.11187  -0.260  0.79519
leisure time      -0.13142    0.12714  -1.034  0.30129
self-esteem        0.45187    0.11105   4.069 4.72e-05
age                0.04298    0.05516   0.779  0.43586
gender             0.19698    0.15494   1.271  0.20360

    Null deviance: 552.4  on 413  degrees of freedom
Residual deviance: 521.4  on 404  degrees of freedom
AIC: 541.4

Number of Fisher Scoring iterations: 4
\end{verbatim}
\end{minipage}
\end{footnotesize}
\end{center}
\end{table}

\begin{table}[htp]
%\label{tab:cancerlm2}
\caption{\label{tab:2}GLM approach: Results from binomial regression of \textit{challenge} on original explanatory composition (in percentages) and further covariates. For explanations see text.}
\begin{center}
\begin{footnotesize}
\begin{minipage}{10cm}
\begin{verbatim}

                   Estimate Std. Error z value Pr(>|z|)
(Intercept)        -0.44332  1.90962  -0.232    0.816
study/work          0.02123  0.01853   1.146    0.252
commuting          -0.00242  0.03916  -0.062    0.951
sleep              -0.00467  0.02077  -0.225    0.822
household duties    0.00073  0.03098   0.024    0.981
leisure time       -0.01892  0.02253  -0.840    0.401
self-esteem         0.44518  0.11046   4.030 5.57e-05
age                 0.03610  0.05418   0.666    0.505
gender              0.16862  0.15198   1.110    0.267

    Null deviance: 552.40  on 413  degrees of freedom
Residual deviance: 524.04  on 405  degrees of freedom
AIC: 542.04

Number of Fisher Scoring iterations: 4
\end{verbatim}
\end{minipage}
\end{footnotesize}
\end{center}
\end{table}

\begin{table}[htp]
%\label{tab:cancerlm2}
\caption{\label{tab:3}Logratio approach: Results from regression with compositional response. Significant regression parameters (at $\alpha=0.05$) marked by asterisk.}
\begin{center}
\begin{footnotesize}
\begin{minipage}{12cm}
\begin{verbatim}

             study/work  commuting     sleep    household  leisure time
(Intercept)  0.60673     -1.17704     1.50068*  -1.27186*    0.96194*
challenge    0.23723*    -0.02216    -0.01922   -0.07174    -0.13237
self-esteem -0.02201      0.02367     0.03083   -0.01959     0.08817*
age          0.04117*    -0.00789     0.00186    0.02723    -0.01801
gender      -0.03217     -0.05116     0.14381*  -0.04412     0.22449*

\end{verbatim}
\end{minipage}
\end{footnotesize}
\end{center}
\end{table}

\begin{table}[htp]
%\label{tab:cancerlm2}
\caption{\label{tab:4}GLM approach: Results from Dirichlet regression with compositional response. Significant regression parameters (at $\alpha=0.05$) marked by asterisk.}
\begin{center}
\begin{footnotesize}
\begin{minipage}{12cm}
\begin{verbatim}

             study/work  commuting     sleep    household  leisure time
(Intercept)  2.05549*     0.94065*    2.59815*   0.98208     2.21778*
challenge    0.08550     -0.05467    -0.06799   -0.08346    -0.14143
self-esteem  0.03210      0.04770     0.06618    0.02641     0.09649*
age          0.00825     -0.01401    -0.01578   -0.00027    -0.02443
gender      -0.03630     -0.03837     0.06755   -0.03792     0.11748*

\end{verbatim}
\end{minipage}
\end{footnotesize}
\end{center}
\end{table}

\begin{table}[htp]
%\label{tab:cancerlm2}
\caption{\label{tab:5}Logratio approach: Results from regression of \textit{leisure time} coordinate on orthogonal coordinates of the explanatory composition and further covariates. For explanations see text.}
\begin{center}
\begin{footnotesize}
\begin{minipage}{10cm}
\begin{verbatim}

                   Estimate Std. Error t value Pr(>|t|)
(Intercept)        0.47988    0.45492   1.055  0.29211
study/work        -0.15646    0.05576  -2.806  0.00526
commuting         -0.20414    0.04377  -4.664 4.22e-06
sleep              0.33976    0.06374   5.330 1.64e-07
household duties   0.05734    0.04304   1.332  0.18351
challenge         -0.09358    0.08937  -1.047  0.29568
self-esteem        0.08353    0.04330   1.929  0.05442
age               -0.01908    0.01991  -0.958  0.33852
gender             0.16760    0.05951   2.817  0.00509

Residual standard error: 0.8544 on 404 degrees of freedom
Multiple R-squared:  0.1619,    Adjusted R-squared:  0.1433
F-statistic: 8.674 on 9 and 404 DF,  p-value: 6.21e-12

\end{verbatim}
\end{minipage}
\end{footnotesize}
\end{center}
\end{table}

\begin{table}[htp]
%\label{tab:cancerlm2}
\caption{\label{tab:6}Standard LM approach: Results from regression of \textit{leisure time} on other compositional parts (in percentages) and further covariates. For explanations see text.}
\begin{center}
\begin{footnotesize}
\begin{minipage}{10cm}
\begin{verbatim}

 Coefficients:
                            Estimate Std. Error t value Pr(>|t|)
(Intercept)                 56.88605    2.96589  19.180  < 2e-16
work/study                  -0.60560    0.02698 -22.446  < 2e-16
commuting                   -0.94166    0.07235 -13.016  < 2e-16
sleep                       -0.51838    0.03738 -13.869  < 2e-16
household duties            -0.63007    0.06094 -10.339  < 2e-16
challenge                   -0.41193    0.48376  -0.852  0.39498
self-esteem                  0.46017    0.23468   1.961  0.05058
age                          0.05062    0.10826   0.468  0.64032
gender                       1.00488    0.31835   3.157  0.00172

Residual standard error: 4.636 on 405 degrees of freedom
Multiple R-squared:  0.6014,    Adjusted R-squared:  0.5935
F-statistic: 76.37 on 8 and 405 DF,  p-value: < 2.2e-16

\end{verbatim}
\end{minipage}
\end{footnotesize}
\end{center}
\end{table}

\vspace*{\fill}\newpage

%%\subsection{\centerline{Figures}}

\begin{figure}
\label{edaplots}
\centerline{\includegraphics[width=\textwidth]{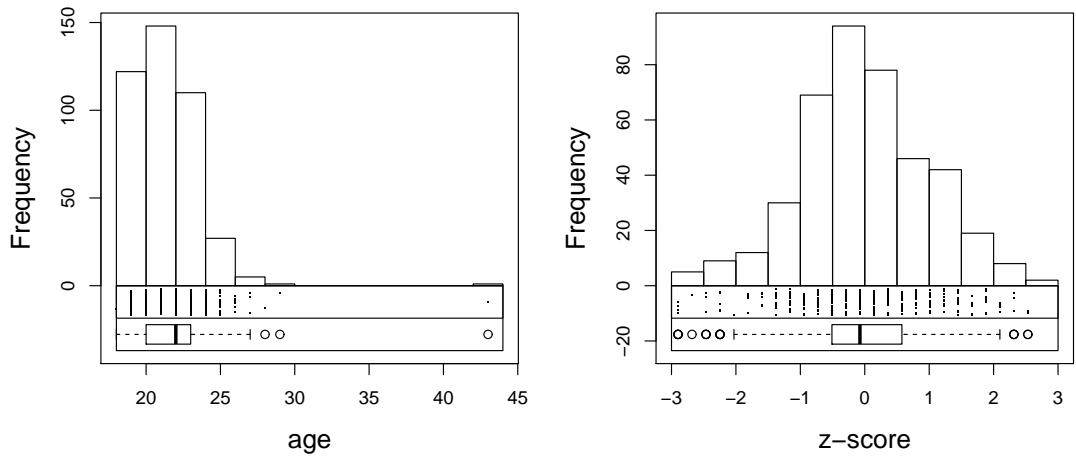}}
\caption{EDA-plots for variables \textit{age} (left) and \textit{self-esteem} (right).}
\end{figure}

%\begin{figure}
%\label{biplots}
%\centerline{\includegraphics[width=\textwidth]{biplots.eps}}
%\caption{Compositional (left) and standard (right) biplots for time budget data. Number codes correspond to single activities (1 $-$ \textit{study/work}, 2 $-$ %\textit{commuting}, 3 $-$ \textit{food}, 4 $-$ \textit{hygiene\&dressing}, 5 $-$ \textit{sleep}, 6 $-$ \textit{household duties}, 7 $-$ \textit{leisure time})}
%\end{figure}

\end{document}